\documentclass[12pt]{amsart}
\usepackage{amssymb}
\usepackage{amsfonts}
\usepackage{amsbsy}
\usepackage{latexsym}
\newtheorem{theorem}{Theorem}
\newtheorem{lemma}[theorem]{Lemma}
\newtheorem{corollary}[theorem]{Corollary}

\theoremstyle{definition}
\newtheorem{definition}[theorem]{Definition}

\begin{document}
\title {Elementary Maps on Triangular algebras}
\author{Xuehan Cheng}
\address{College of  Mathematics and Information, Ludong University,  Yantai, 264025 P. R. China}

\author{Wu Jing}
\address{Department of Mathematics and Computer Science, Fayetteville State University,  Fayetteville, NC 28301}
 \email{wjing@uncfsu.edu}

\thanks{This work is partially supported by NNSF of China (No. 10671086) and NSF of Ludong University (No. LY20062704).}
\date{June 6, 2007}
\begin{abstract}
In this note we prove that  elementary surjective maps on triangular algebras are automatically additive. \end{abstract}
\maketitle

The study of elementary maps was initiated by Bre$\check{\textrm {s}}$ar and $\check{\textrm{S}}$erml. Following  (\cite
{br1115}), elementary maps are defined as follows.
\begin{definition}
Let $\mathcal {R}$ and $\mathcal {R}^{\prime }$ be two rings. Suppose that $M\colon \mathcal {R}\rightarrow \mathcal {R}^{\prime}$ and
$M^*\colon \mathcal {R}^{\prime}\rightarrow \mathcal {R}$ are two
maps.  Call the ordered pair $(M, M^*)$ an \textit {elementary
map} of $\mathcal {R}\times \mathcal {R}^{\prime}$ if
\begin{displaymath}
 \left\{ \begin{array}{ll}
 M(aM^*(x)b)=M(a)xM(b),\\
 M^*(xM(a)y)=M^*(x)aM^*(y)
 \end{array}\right.
\end{displaymath}
 for all $a, b\in \mathcal {R}$ and $x, y\in \mathcal {R}^{\prime}$.
 \end{definition}

 Notice that no additivity of $M$ and $M^*$ is required in the above definition.

 It is very interesting that elementary maps on some rings as well as operator algebras are automatically additive. In recent years, the study of additivity of maps on rings and operator algebras has become an active topic which has enriched the theory of maps on rings and operator algebras.

The first result about the additivity of maps on rings was due to  Martindale III. In \cite {ma695}, he proved the following result.
\begin{theorem}(\cite{ma695}) Let $\mathcal {R}$ be a ring containing a family $\{ e_{\alpha }:\alpha \in \Lambda \} $ of idempotents which satisfies:

(i) $x\mathcal {R}=\{ 0\} $ implies $x=0$;

(ii) If $e_{\alpha }\mathcal {R}x=\{ 0\} $ for each $\alpha \in \Lambda $, then $x=0$ (and hence $\mathcal {R}x=\{ 0\} $ implies $x=0$);

(iii) For each $\alpha \in \Lambda $, $e_{\alpha }xe_{\alpha }\mathcal {R}(1-e_{\alpha })=\{ 0\} $ implies $e_{\alpha }xe_{\alpha }=0$.

Then any multiplicative bijective map from $\mathcal {R}$ onto an arbitrary ring $\mathcal {R}^{\prime }$ is additive.
\end{theorem}

As a corollary, every multiplicative bijective map from a prime ring containing a nontrivial idempotent onto an arbitrary ring is necessarily additive.

It should be mentioned here that the proof of \cite {ma695} has become a quite standard argument and been applied widely in dealing with the additivity of a large number of maps on rings and operator algebras (see \cite {ji}-\cite{lu2273}).

 The aim of this note is   to study the additivity  of elementary maps on triangular algebras. We will show that every elementary map on triangular operator is additive. Note that, different from \cite {ma695},  we do not require the existence of  nontrivial idempotents. However, we still follow the line of \cite {ma695}.

 Recall that a \textit{triangular algebra} $Tri(\mathcal {A}, \mathcal {M}, \mathcal {B})$ is  an algebra of the form

 \begin{displaymath}Tri(\mathcal {A}, \mathcal {M}, \mathcal {B})=\{
\left( \begin{array}{ll}
a & m\\
0 & b
\end{array}\right): a\in \mathcal {A}, m\in \mathcal {M}, b\in \mathcal {B}\}
\end{displaymath}
under the usual matrix operations, where $\mathcal {A}$ and $\mathcal {B}$ are two algebras over a commutative ring $\mathcal {R}$, and $\mathcal {M}$ is an $(\mathcal {A}, \mathcal {B})$-bimodule which is faithful as a left $\mathcal {A}$-module and also as a right $\mathcal {B}$-module (see \cite{cheung}).

We are ready to state our result of this note.
\begin{theorem} \label{theorem} Let $\mathcal {R}^{\prime }$ be an arbitrary ring. Let $\mathcal {A}$ and $\mathcal {B}$ be two algebras over a commutative ring $\mathcal {R}$, $\mathcal {M}$ a faithful $(\mathcal {A}, \mathcal {B})$-bimodule, and $\mathcal {T}$ be the triangular algebra $Tri(\mathcal {A}, \mathcal {M}, \mathcal {B})$. Suppose that algebras $\mathcal {A}$ and $\mathcal {B}$ satisfy:

(i) If $a\mathcal {A}=\{ 0\} $, or $\mathcal {A}a=\{ 0\} $, then $a=0$;

(ii) If $b\mathcal {B}=\{ 0\} $, or $\mathcal {B}b=\{ 0\} $, then $b=0$.

Suppose that $(M, M^*)$ is an elementary map on $\mathcal {T}\times \mathcal {R}^{\prime }$, and both $M$ and $M^*$ are surjective. Then both $M$ and $M^*$ are additive.
\end{theorem}

 For the sake of clarity, we divide the proof into a series of lemmas. We begin with the following trivial one.

\begin{lemma} $M(0)=0$ and $M^*(0)=0$.
\end{lemma}
\begin{proof} We have
$M(0)=M(0M^*(0)0)=M(0)0M(0)=0$.

Similarly, $M^*(0)=M^*(0M(0)0)=M^*(0)0M^*(0)=0$.
\end{proof}

In what follows, we set \begin{displaymath} \mathcal {T}_{11}=\{
\left( \begin{array}{ll}
a & 0\\
0 & 0
\end{array}\right): a\in \mathcal {A} \},
\end{displaymath}
\begin{displaymath}\mathcal {T}_{12}=\{
\left( \begin{array}{ll}
0 & m\\
0 & 0
\end{array}\right): m\in \mathcal {M} \},
\end{displaymath}
and
\begin{displaymath} \mathcal {T}_{22}=\{
\left( \begin{array}{ll}
0 & 0\\
0 & b
\end{array}\right): b\in \mathcal {B} \}.
\end{displaymath}
Then we may write $\mathcal {T}=\mathcal {T}_{11}\oplus \mathcal {T}_{12}\oplus \mathcal {T}_{22}$, and every element $a\in\mathcal {T}$ can be written as $a=a_{11}+a_{12}+a_{22}$. Note that notation $a_{ij}$ denotes an arbitrary element of
 $\mathcal {T}_{ij}$.

The following result shows that both $M$ and $M^*$ are bijective.
  \begin{lemma}
Both $M$ and $M^*$ are injective.
\end{lemma}
\begin{proof}   Suppose that $M(a)=M(b)$ for some $a$ and $b$ in $\mathcal {T}$.
We write  $a=a_{11}+a_{12}+a_{22}$ and $b=b_{11}+b_{12}+b_{22}$.

For arbitrary $x$ and $y$ in $\mathcal {R}^{\prime }$, we have
$$M^*(x)aM^*(y)=M^*(xM(a)y)=M^*(aM(b)y)=M^*(x)bM^*(y).$$
This, by the surjectivity of $M^*$, is equivalent to \begin{equation}\label{01}
sat=sbt\end{equation}
for arbitrary $s, t\in \mathcal {T}$.

In particular, letting $s=s_{11}, t=t_{ii}\in \mathcal {T}_{11}$ in equality (\ref{01}), we get $s_{11}a_{11}t_{11}=s_{11}b_{11}t_{11}$. And so,  by condition (i) in Theorem \ref{theorem}, $a_{11}=b_{11}$.

Similarly, we can get $a_{22}=b_{22}$ by letting $s=s_{22}$ and $t=t_{22}$ in identity (\ref{01}).

We now show that $a_{12}=b_{12}$. Setting $s=s_{11}\in \mathcal {T}_{11}$ and $t=t_{22}\in \mathcal {T}_{22}$, then equality (\ref{01}) becomes $s_{11}at_{22}=s_{11}bt_{22}$, i.e., $s_{11}a_{12}t_{22}=s_{11}b_{12}t_{22}$. Therefore $a_{12}=b_{12}$ as $\mathcal {M}$ is a faithful $(\mathcal {A}, \mathcal {B})$-bimodule.

To complete the proof, it remains to show that $M^*$ is injective. Let $x$ and $y$ be in $\mathcal {R}^{\prime }$ such that $M^*(x)=M^*(y)$. Now for any $a, b\in \mathcal {T}$, we have
\begin{eqnarray*}
& & M^*M(a)M^{-1}(x)M^*M(b)\\
&=&M^*(M(a)MM^{-1}(x)M(b))\\
&=&M^*(M(a)xM(b))\\
&=&M^*M(aM^*(x)b)\\
&=&M^*M(aM^*(y)b)\\
&=&M^*(M(a)yM(b))\\
&=&M^*(M(a)MM^{-1}(y)M(b))\\
&=&M^*M(a)M^{-1}(y)M^*(b).
\end{eqnarray*}
Thus
$$M^*M(a)M^{-1}(x)M^*M(b)=M^*M(a)M^{-1}(y)M^*M(b).$$
Equivalently,
$$sM^{-1}(x)t=sM^{-1}(y)t$$
for any $s, t\in \mathcal {T}$ since $M^*M$ is surjective.

It follows from the same argument above that $M^{-1}(x)=M^{-1}(y)$, and so $x=y$, as desired.
\end{proof}

 \begin{lemma}\label{inverse}
The pair $(M^{*^{-1}}, M^{-1})$ is a Jordan elementary map on
$\mathcal {T}\times \mathcal {R}^{\prime}$. That is,
\begin{displaymath}
 \left\{ \begin{array}{ll}
 M^{*^{-1}}(aM^{-1}(x)b)=M^{*^{-1}}(a)xM^{*^{-1}}(b),\\
 M^{-1}(xM^{*^{-1}}(a)y)=M^{-1}(x)aM^{-1}(y)
\end{array}\right.
\end{displaymath}
for all $a, b\in \mathcal {T}$ and  $x, y\in \mathcal {R}^{\prime }$.
\end{lemma}
\begin{proof}
The first identity follows from the following observation.
\begin{eqnarray*}
& &M^*(M^{*^{-1}}(a)xM^{*^{-1}}(b))\\
&=&M^*(M^{*^{-1}}(a)MM^{-1}(x)M^{*^{-1}}(b))\\
&=&aM^{-1}(x)b.
\end{eqnarray*}
  The second one goes similarly.
\end{proof}

The following result will be used frequently in this note.
\begin{lemma}\label{add}
Let $a, b, c\in \mathcal {R}$ such that $M(c)=M(a)+M(b)$. Then for any $s, t\in \mathcal {T}$ we have
$$M^{*^{-1}}(sct)=M^{*^{-1}}(sat)+M^{*^{-1}}(sbt).$$
 \end{lemma}
\begin{proof}  By Lemma \ref{inverse}, we have
\begin{eqnarray*}
M^{*^{-1}}(sct)&=&M^{*^{-1}}(sM^{-1}M(c)t)\\
&=&M^{*^{-1}}(s)M(c)M^{*^{-1}}(t)\\
&=&M^{*^{-1}}(s)(M(a)+M(b))M^{*^{-1}}(t)\\
&=&(M^{*^{-1}}(s)M(a)M^{*^{-1}}(t))+(M^{*^{-1}}(s)M(b)M^{*^{-1}}(t))\\
&=&M^{*^{-1}}(sat)+M^{*^{-1}}(sbt).
\end{eqnarray*}
\end{proof}

\begin{lemma}\label{lemma1112}
Let $a_{11}\in \mathcal {T}_{11}$ and $b_{12}\in \mathcal {T}_{12}$,
then

(i) $M(a_{11}+b_{12})=M(a_{11})+M(b_{12})$;

(ii)
$M^{*^{-1}}(a_{11}+b_{12})=M^{*^{-1}}(a_{11})+M^{*^{-1}}(b_{12})$.
\end{lemma}
\begin{proof} We only prove (i).
We choose  $c\in \mathcal {T}$   such that $M(c)=M(a_{11})+M(b_{12})$. For arbitrary $s_{11}\in \mathcal {R}_{11}$ and $t_{22}\in \mathcal {T}_{22}$, by Lemma \ref{add}, we have
$$
M^{*^{-1}}(s_{11}ct_{22})=M^{*^{-1}}(s_{11}a_{11}t_{22})+M^{*^{-1}}(s_{11}b_{12}t_{22})=M^{*^{-1}}(s_{11}b_{12}t_{22}).$$
It follows that $s_{11}ct_{22}=s_{11}b_{12}t_{22}$, i.e., $s_{11}c_{12}t_{22}=s_{11}b_{12}t_{22}$, which yields that $c_{12}=b_{12}$.

Now for any $s_{11}$ and $t_{11}$ in $\mathcal {T}_{11}$, we have
$$M^{*^{-1}}(s_{11}ct_{11})=M^{*^{-1}}(s_{11}a_{11}t_{11})+M^{*^{-1}}(s_{11}b_{12}t_{11})=M^{*^{-1}}(s_{11}a_{11}t_{11}).$$  This implies that $c_{11}=a_{11}$.

Similarly, for any $s_{22}, t_{22}\in \mathcal {T}_{22}$, we obtain
 $$M^{*^{-1}}(s_{22}ct_{22})=M^{*^{-1}}(s_{22}a_{11}t_{22})+M^{*^{-1}}(s_{22}b_{12}t_{22})=0.$$
 Hence $c_{22}=0$ follows from the fact that $s_{22}ct_{22}=0$..
\end{proof}

Similarly, we can get the following result.
\begin{lemma}\label{lemma2212}
Let $a_{22}\in \mathcal {T}_{22}$ and $b_{12}\in \mathcal {T}_{12}$,
 then

(i) $M(a_{22}+b_{12})=M(a_{22})+M(b_{12})$;

(ii)
$M^{*^{-1}}(a_{22}+b_{12})=M^{*^{-1}}(a_{22})+M^{*^{-1}}(b_{12})$.
\end{lemma}

\begin{lemma}\label{lemma111112}
For any $t_{11}, a_{11}\in \mathcal {T}_{11}$, $b_{12}, c_{12}\in \mathcal {T}_{12}$, and $d_{22}\in \mathcal {T}_{22}$, we have

(i) $M(t_{11}a_{11}b_{12}+t_{11}c_{12}d_{22})=M(t_{11}a_{11}b_{12})+M(t_{11}c_{12}d_{22})$;

(ii) $M^{*^{-1}}(t_{11}a_{11}b_{12}+t_{11}c_{12}d_{22})=M^{*^{-1}}(t_{11}a_{11}b_{12})+M^{*^{-1}}(t_{11}c_{12}d_{22}).$
\end{lemma}
\begin{proof} We only prove (i).
Using Lemma \ref{lemma1112} and Lemma \ref{lemma2212}, we compute
\begin{eqnarray*}
& &M(t_{11}a_{11}b_{12}+t_{11}c_{12}d_{22})\\
&=&M(t_{11}(a_{11}+c_{12})(b_{12}+d_{22}))\\
&=&M(t_{11}M^*M^{*^{-1}}(a_{11}+c_{12})(b_{12}+d_{22}))\\
&=&M(t_{11})M^{*^{-1}}(a_{11}+c_{12})M(b_{12}+d_{22})\\
&=&M(t_{11})M^{*^{-1}}(a_{11})M(b_{12})+M(t_{11})M^{*^{-1}}(a_{11})M(d_{22})\\
& &+M(t_{11})M^{*^{-1}}(c_{12})M(b_{12})+M(t_{11})M^{*^{-1}}(c_{12})M(d_{22})\\
&=&M(t_{11})(M^{*^{-1}}(a_{11})+M^{*^{-1}}(c_{12}))M(b_{12})+M(t_{11})(M^{*^{-1}}(a_{11})+M^{*^{-1}}(c_{12}))M(d_{22})\\
&=&M(t_{11})M^{*^{-1}}(a_{11}+c_{12})M(b_{12})+M(t_{11})M^{*^{-1}}(a_{11}+c_{12})M(d_{22})\\
&=&M(t_{11}(a_{11}+c_{12})b_{12})+M(t_{11}(a_{11}+c_{12})d_{22})\\
&=&M(t_{11}a_{11}b_{12})+M(t_{11}c_{12}d_{22}).
\end{eqnarray*}
\end{proof}

\begin{lemma}\label{lemma12}
Both $M$ and $M^{*^{-1}}$ are additive on $\mathcal {T}_{12}$.
\end{lemma}
\begin{proof}
Let  $a_{12}$ and $b_{12}$ be in $\mathcal {T}_{12}$. We pick $c\in \mathcal {T}$ such that $M(c)=M(a_{12})+M(b_{12})$.

For arbitrary $t_{11}, s_{11}\in \mathcal {T}_{11}$, by Lemma \ref{add}, we have
$$M(t_{11}cs_{11})=M(t_{11}a_{12}s_{11})+M(t_{11}b_{12}s_{11})=0,$$
this implies that $t_{11}cs_{11}=0$, and so $c_{11}=0$.

Similarly, we can get $c_{22}=0$.

We now show that $c_{12}=a_{12}+b_{12}$. For any $t_{11}, r_{11}\in \mathcal {T}_{11}$ and $s_{22}\in \mathcal {T}_{22}$, by Lemma \ref{add} and Lemma \ref{lemma111112}, we obtain
\begin{eqnarray*}
M(r_{11}t_{11}cs_{22})&=&M(r_{11}t_{11}a_{12}s_{22})+M(r_{11}t_{11}a_{12}s_{22})\\
&=&M(r_{11}t_{11}a_{12}s_{22}+r_{11}t_{11}b_{12}s_{22})\\
&=&M(r_{11}t_{11}(a_{12}+b_{12})s_{22}).
\end{eqnarray*}
It follows that $$r_{11}t_{11}cs_{22}=r_{11}t_{11}(a_{12}+b_{12})s_{22}.$$
Equivalently, $$r_{11}t_{11}c_{12}s_{22}=r_{11}t_{11}(a_{12}+b_{12})s_{22}.$$
Then we get $c_{12}=a_{12}+b_{12}$.

With the similar argument, one can see that $M^{*^{-1}}$ is also additive on $\mathcal {T}_{12}$.
\end{proof}
\begin{lemma}\label{lemma11}
Both $M$ and $M^{*^{-1}}$ are additive on $\mathcal {T}_{11}$.
\end{lemma}
\begin{proof} We only show the additivity of $M$ on $\mathcal {T}_{11}$.
Suppose that $a_{11}$ and $b_{11}$ are two elements of $\mathcal {T}_{11}$. Let $c\in \mathcal {T}$ be chosen satisfying $M(c)=M(a_{11})+M(b_{11})$. Now for any $t_{22}, s_{22}\in \mathcal {T}_{22}$, we have $$M(t_{22}cs_{22})=M(t_{22}a_{11}s_{22})+M(t_{22}b_{11}s_{22})=0.$$
Consequently, $t_{22}cs_{22}=0$, i.e., $t_{22}c_{22}s_{22}=0$, and so $c_{22}=0$.

Similarly, we can infer that $c_{12}=0$.

To complete the proof, we need to show that $c_{11}=a_{11}+b_{11}$. For each $t_{11}\in \mathcal {T}_{11}$ and $s_{12}\in \mathcal {T}_{12}$, we consider

$$M(t_{11}cs_{12})=M(t_{11}a_{11}s_{12})+M(t_{11}b_{11}s_{12})=M(t_{11}a_{11}s_{12}+t_{11}b_{11}s_{12}).$$
Note that in the last equality we apply Lemma \ref{lemma12}. It follows that $$t_{11}cs_{12}=t_{11}a_{11}s_{12}+t_{11}b_{11}s_{12}.$$
This leads to $c_{11}=a_{11}+b_{11}$, as desired.
\end{proof}
\begin{lemma}\label{lemma22}
$M$ and $M^{*^{-1}}$ are additive on $\mathcal {T}_{22}$.
\end{lemma}
\begin{proof}
Suppose that $a_{22}$ and $b_{22}$ are in $\mathcal {T}_{22}$. We choose $c\in \mathcal {T}$ such that $M(c)=M(a_{22})+M(b_{22})$. For any $t_{12}\in\mathcal {T}_{12}$ and $s_{22}\in \mathcal {T}_{22}$, using Lemma \ref{lemma12}, we have
$$M(t_{12}cs_{22})=M(t_{12}a_{22}s_{22})+M(t_{12}b_{22}s_{22})=M(t_{12}(a_{22}+b_{22})s_{22}).$$
Accordingly, $t_{12}cs_{22}=t_{12}(a_{22}+b_{22})s_{22}$, which yields that $c_{22}=a_{22}+b_{22}$.

With the similar argument, we can verify that $c_{11}=c_{12}=0$.

The additivity of $M^{*^{-1}}$ on $\mathcal {T}_{22}$ follows  easily.
\end{proof}
\begin{lemma}\label{lemma111222}
For any $a_{11}\in \mathcal {T}_{11}$, $b_{12}\in \mathcal {T}_{12}$, and $c_{22}\in\mathcal {T}_{22}$, the following are true.

(i) $M(a_{11}+b_{12}+c_{22})=M(a_{11})+M(b_{12})+M(c_{22})$;

(ii) $M^{*^{-1}}(a_{11}+b_{12}+c_{22})=M^{*^{-1}}(a_{11})+M^{*^{-1}}(b_{12})+M^{*^{-1}}(c_{22}).$

\end{lemma}

\begin{proof} We only prove (i).
Let $d\in \mathcal {T}$ be an element satisfying $M(d)=M(a_{11})+M(b_{12})+M(c_{22})$.  For any $s, t\in \mathcal {T}$, using Lemma \ref{add} twice, we can arrive at
\begin{equation}\label{22}
M(sdt)=M(sa_{11}t)+M(sb_{12}t)+M(sc_{22}t).
\end{equation}

Letting $s=s_{11}$ and $t=t_{11}$ in the above equality, we  get $d_{11}=a_{11}$.

In the same fashion for $s=s_{22}$ and $t=t_{22}$ in equality (\ref{22}), we can infer that $d_{22}=c_{22}$.

Finally, considering  $s=s_{11}$ and $t=t_{22}$ in equality (\ref{22}), we see that $d_{12}=b_{12}$. Thus, $d=a_{11}+b_{12}+c_{22}$, which completes the proof.
\end{proof}
 We now prove our main result.

  \noindent \textbf{Proof of Theorem \ref{theorem}} We first show the additivity of  $M$. Let $a=a_{11}+a_{12}+a_{22}$ and $b=b_{11}+b_{12}+b_{22}$ be two arbitrary elements of $\mathcal {T}$. We have
 \begin{eqnarray*}
 & &M(a+b)\\
 &=&M((a_{11}+b_{11})+(a_{12}+b_{12})+(a_{22}+b_{22}))\\
 &=&M(a_{11}+b_{11})+M(a_{12}+b_{12})+M(a_{22}+b_{22})\\
 &=&M(a_{11})+M(b_{11})+M(a_{12})+M(b_{12})+M(a_{22})+M(b_{22})\\
 &=&(M(a_{11})+M(a_{12})+M(a_{22}))+(M(b_{11})+M(b_{12})+M(b_{22}))\\
 &=&M(a_{11}+a_{12}+a_{22})+M(b_{11}+b_{12}+b_{22})\\
 &=&M(a)+M(b).
 \end{eqnarray*}
 That is, $M$ is additive.

 We now turn to prove that $M^*$ is additive. For any $x, y\in \mathcal {R}^{\prime}$, there exist $c=c_{11}+c_{12}+c_{22}$ and $d=d_{11}+d_{12}+d_{22}$ in $\mathcal {R}$ such that $c=M^*(x+y)$ and $d=M^*(x)+M^*(y)$.

 For arbitrary $s, t\in \mathcal {T}$, by the additivity of $M$, we compute
 \begin{eqnarray*}
  M(sct)&=&M(sM^*(x+y)t)\\
  &=&M(s)(x+y)M(t)\\
  &=&M(s)xM(t)+M(s)yM(t)\\
  &=&M(sM^*(x)t)+M(sM^*(y)t)\\
  &=&M(sM^*(x)t+sM^*(y)t)\\
  &=&M(s(M^*(x)+M^*(y))t)\\
  &=&M(sdt),
 \end{eqnarray*}
 which implies that $sct=sdt$. Furthermore, we get $c=d$, i.e., $M^*(x+y)=M^*(x)+M^*(y)$.

In particular, if both $\mathcal {A}$ and $\mathcal {B}$ are unital algebras, we have
\begin{corollary}
Let $\mathcal {R}^{\prime }$ be an arbitrary ring. Let $\mathcal {A}$ and $\mathcal {B}$ be two unital algebras over a commutative ring $\mathcal {R}$, $\mathcal {M}$ a faithful $(\mathcal {A}, \mathcal {B})$-bimodule, and $\mathcal {T}$ be the triangular algebra $Tri(\mathcal {A}, \mathcal {M}, \mathcal {B})$.
Suppose that $(M, M^*)$ is an elementary map on $\mathcal {T}\times \mathcal {R}^{\prime }$, and both $M$ and $M^*$ are surjective. Then both $M$ and $M^*$ are additive.\end{corollary}

Recall that a subalgebra of $B(X)$ is called a \textit {standard operator} if it contains all finite rank operators, where $B(X)$ is the algebra of all bounded linear operator on a Banach space $X$.

 We complete this note by considering   elementary maps on triangular algebras provided $\mathcal {A}$ and $\mathcal {B}$ are standard operator algebras.

\begin{corollary}
Let $\mathcal {R}^{\prime }$ be an arbitrary ring. Let $\mathcal {A}$ and $\mathcal {B}$ be two standard operator algebras over a Banach space $X$, $\mathcal {M}$ a faithful $(\mathcal {A}, \mathcal {B})$-bimodule, and $\mathcal {T}$ be the triangular algebra $Tri(\mathcal {A}, \mathcal {M}, \mathcal {B})$.
Suppose that $(M, M^*)$ is an elementary map on $\mathcal {T}\times \mathcal {R}^{\prime }$, and both $M$ and $M^*$ are surjective. Then both $M$ and $M^*$ are additive.\end{corollary}
 \bibliographystyle{amsplain}

\begin{thebibliography}{10}

\bibitem {br1115} M. Bre$\check{\textrm {s}}$ar, P.  $\check{\textrm{S}}$erml, Elementary operators, \textit {Proc. Roy. Soc. Edinburgh Ser A.}, \textbf {129}, (1999), 1115--1135.


\bibitem{cheung} W. S. Cheung, Commuting maps on triangular algebras, \textit {J. London Math. Soc.}, \textbf {63} (2001), 117--127.

 \bibitem {ji} P. Ji,   Jordan maps on triangular algebras, \textit {Linear Algebra Appl.}, (to appear).

\bibitem {jing1} W. Jing, Aditivity of Jordan elementary maps on rings, preprint (2007).


\bibitem {jing2} W. Jing, Jordan triple elementary maps on rings, preprint (2007).

\bibitem {li237} P. Li, W. Jing, Jordan elementary maps on rings, \textit {Linear Algebra Appl.}, \textbf {382} (2004), 237--245.

\bibitem{lilu} P. Li, F. Lu, Additivity of elementary maps on rings, \textit{Comm. Algebra}, \textbf{32} (2004), 3725--3737.

\bibitem{lu123} F. Lu, Additivity of Jordan maps on standard operator algebras, \textit {Linear Algebra Appl.}, \textbf{357} (2002), 123--131.


\bibitem{lu2273} F. Lu, Jordan maps on associative algebras, \textit {Comm. Algebra}, \textbf {31} (2003), 2273--2286.
\bibitem{ma695} W. S. Martindale III, When are multiplicative mappings additive?\textit {Proc. Amer. Math. Soc.}, \textbf {21} (1969) 695--698.
 \end{thebibliography}

\end{document}